\documentclass[a4paper,10pt,twoside]{article}

\usepackage[utf8]{inputenc}
\usepackage{amsmath, amssymb, amsthm}
\usepackage{amstext, amsfonts, a4}
\usepackage{dsfont}
\usepackage{latexsym}
\usepackage{mathrsfs}
\usepackage[all,cmtip]{xy}
\usepackage{graphicx}
\usepackage{enumerate}
\usepackage{hyperref,color}
\usepackage{stmaryrd}
\usepackage[shortlabels]{enumitem}

\usepackage{url}
\usepackage[pdftex,a4paper,left=3cm,right=3cm,top=3cm,bottom=3cm]{geometry}

\def\sV{\mathscr{V}}

\def\bbC{\mathbb{C}}
\def\bbF{\mathbb{F}}

\def\bbQ{\mathbb{Q}}

\def\bbZ{\mathbb{Z}}

\def\cI{\mathcal{I}}

\DeclareMathOperator{\End}{End}

\DeclareMathOperator{\ind}{ind}

\newtheorem{counter}[subsection]{$\!\!$}

\newenvironment{Lem}{\begin{counter} {\bf Lemma.}}{\end{counter}}
\newenvironment{Cor}{\begin{counter} {\bf Corollary.}}{\end{counter}}
\newenvironment{Th}{\begin{counter} {\bf Theorem.}}{\end{counter}}

\newtheorem{counter*}[subsubsection]{$\!\!$}
\newenvironment{Def*}{\begin{counter*} {\bf Definition.}}{\end{counter*}}
\newenvironment{Not*}{\begin{counter*} \rm {\bf Notation.}}{\end{counter*}}
\newenvironment{Notss*}{\begin{counter*} \rm {\bf Notations.}}{\end{counter*}}
\newenvironment{DefNot*}{\begin{counter*} \rm {\bf Definition-Notation.}}{\end{counter*}}
\newenvironment{Nots*}{\begin{counter*} \rm {\bf Notations.}}{\end{counter*}}
\newenvironment{Prop*}{\begin{counter*} {\bf Proposition.}}{\end{counter*}}
\newenvironment{Lem*}{\begin{counter*} {\bf Lemma.}}{\end{counter*}}
\newenvironment{Cor*}{\begin{counter*} {\bf Corollary.}}{\end{counter*}}
\newenvironment{Th*}{\begin{counter*} {\bf Theorem.}}{\end{counter*}}
\newenvironment{Rem*}{\begin{counter*} \rm {\bf Remark.}}{\end{counter*}}
\newenvironment{Ex*}{\begin{counter*} \rm {\bf Example.}}{\end{counter*}}
\newenvironment{Exs*}{\begin{counter*} \rm {\bf Examples.}}{\end{counter*}}
\newenvironment{Pt*}{\begin{counter*} \rm}{\end{counter*}}
\newenvironment{Q*}{\begin{counter*} \rm {\bf Question.}}{\end{counter*}}

\title{\textbf{\huge{Irreducible mod $p$ Lubin-Tate $(\varphi,\Gamma)$-modules}}}

\author{Cédric Pépin and Tobias Schmidt}
\date{\today}

\begin{document}

\maketitle

\begin{abstract}
Let $F$ be a finite extension of $\bbQ_p$. We determine the Lubin-Tate 
$(\varphi,\Gamma)$-modules associated to
the absolutely irreducible mod $p$ representations of the absolute Galois group ${\rm Gal}(\overline{F}/F)$.
\end{abstract}

\tableofcontents

\section{Introduction}

Let $F$ be a finite extension of $\bbQ_p$ with ring of integers $o_F$, residue field $\bbF_q$ and uniformizer $\pi\in o_F$. Let $\overline{F}$ be an algebraic closure of $F$ and let
${\rm Gal}(\overline{F}/F)$ be the absolute Galois group of $F$. 
By adapting the well-known formalism of Fontaine for the cyclotomic case, Kisin-Ren explained in \cite{KR09} (see also the detailed exposition by Schneider \cite{Sch17}) how to build an equivalence between the category 
of continuous representations of ${\rm Gal}(\overline{F}/F)$ in finitely generated $o_F$-modules and a category of \' etale Lubin-Tate $(\varphi,\Gamma)$-modules. Let $k/\bbF_q$ be a finite extension. Via reduction modulo $\pi$ and extension of scalars, one deduces an equivalence of categories 
between smooth representations of ${\rm Gal}(\overline{F}/F)$
in finite dimensional $k$-vector spaces and a category of Lubin-Tate 
$(\varphi,\Gamma)$-modules over the Laurent series ring $k((t))$. 
When $F=\bbQ_p$ and in the cyclotomic case, the $(\varphi,\Gamma)$-modules
corresponding to the $n$-dimensional {\it absolutely irreducible} mod $p$ Galois representations have been explicitly calculated by Berger \cite{Be10} and then used by him, in the case of $n=2$, to give a direct proof of the compatibility of Colmez' $p$-adic local Langlands correspondence with Breuil's mod $p$ correspondence for the group ${\rm GL}_2(\bbQ_p)$ in the irreducible case.
In view of extending such results to more general base fields $F\neq \bbQ_p$, we propose in this note to explicitly calculate the Lubin-Tate $(\varphi,\Gamma)$-modules corresponding to the absolutely irreducible mod $p$ representations of ${\rm Gal}(\overline{F}/F)$ for $F\neq \bbQ_p$, thereby generalizing Berger's result. As a method of proof, we adapt Berger's strategy to the Lubin-Tate setting. 

\vskip5pt

In \cite{GK18} (generalizing \cite{GK16} for $F=\bbQ_p$) Grosse-Kl\"onne constructs a fully faithful exact 
functor from the category of so-called supersingular modules
for the pro-$p$ Iwahori-Hecke algebra over $k$ of the group ${\rm GL}_n(F)$ to the category of Lubin-Tate $(\varphi,\Gamma)$-modules over $k((t))$. 
It induces a bijection between the absolutely irreducible objects of rang $n$ on both sides. In \cite{PS1} we show, as an application of the results in this note, how to geometrically construct an inverse map to Grosse-Kl\"onne's bijection in the case $n=2$.

\vskip5pt 

The second author thanks L. Berger for answering some questions on $(\varphi,\Gamma)$-modules.

\section{Galois representations and Lubin-Tate $(\varphi,\Gamma)$-modules}

Let $F_n/F$ be the unramified extension of degree $n$. The irreducible smooth $\overline{\bbF}_q$-representations of ${\rm Gal}(\overline{F}/F) $ of dimension $n$ are 
given by the representations 

$$ \ind^{ {\rm Gal}(\overline{F}/F) }_{ {\rm Gal}(\overline{F}/F_n) }(\chi)$$
smoothly induced from the regular $\overline{\bbF}_q$-characters $\chi$ of ${\rm Gal}(\overline{F}/F_n).$ The ${\rm Gal}(F_n/F)$-conjugates of $\chi$ induce isomorphic representations and there are no other isomorphisms between the representations \cite{V94}.

\vskip5pt 
Let $\pi\in o_F$ be a uniformizer and let $q=p^f$. Let $\pi_{nf}\in \overline{F}$ be an element such that $\pi_{nf}^{q^n-1}=-\pi$. We then have 
Serre's fundamental character of level $nf$
$$\omega_{nf}: {\rm Gal}(\overline{F}/F_{n}) \longrightarrow \bbF_{q^n}^{\times}$$ given by 
$ g\mapsto g(\pi_{nf})/\pi_{nf} \in \mu_{q^n-1}(\overline{F})$ followed by reduction mod $\pi$, cf. \cite{Se72}.
One has $\omega_{nf}^{\frac{q^n-1}{q-1}}=\omega_f|_{{\rm Gal}(\overline{F}/F_{n})}$.

\vskip5pt 

Let $\cI\subset {\rm Gal}(\overline{F}/F)$ be the inertia subgroup and choose an element $\varphi\in {\rm Gal}(\overline{F}/F)$ lifting the Frobenius $x\mapsto x^q$ on $ {\rm Gal}(\overline{\bbF}_q/\bbF_q)$. Since $\omega_{f}: \cI \rightarrow\bbF^{\times}_q$ is surjective \cite[Prop. 2]{Se72}, we may and will assume $\omega_{f}(\varphi)=1$.

\vskip5pt 

A character $\omega_{nf}^h$ for $1\leq h \leq q^n-2$ is regular if and only if its conjugates $\omega_{nf}^{h}, \omega_{nf}^{qh},..., \omega_{nf}^{q^{n-1}h}$ are all distinct. Equivalently, if and only if $h$ is $q$-primitive, that is, there is no $d<n$
such that $h$ is a multiple of $(q^n-1)/(q^d-1)$.
 The representation 
 $\ind^{{\rm Gal}(\overline{F}/F)}_{ {\rm Gal}(\overline{F}/F_n) }(\omega_{nf}^h)$ is then defined over $\bbF_{q^n}$. It has a basis 
 $\{ v_0,...,v_{n-1}\}$ of eigenvectors for the characters 
 $\omega_{nf}^{h}, \omega_{nf}^{qh},..., \omega_{nf}^{q^{n-1}h}$ 
 of ${\rm Gal}(\overline{F}/F_n)$ such that $\varphi(e_i)=e_{i-1}$ and $\varphi(e_0)=e_{n-1}$. In particular, its determinant coincides with 
 $\omega^h_f$ on the subgroup
 ${\rm Gal}(\overline{F}/F_{n})$ and takes $\varphi$ to $(-1)^{n-1}$.
  
For $\lambda\in\overline{\bbF}_q^{\times}$, let $\mu_{\lambda}$ be the unramified character of ${\rm Gal}(\overline{F}/F)$ sending $\varphi$ to 
$\lambda^{-1}$. Fix $\delta$ with $\delta^n=(-1)^{n-1}$.  The representation
 $$ \ind (\omega_{nf}^h):= ( \ind^{{\rm Gal}(\overline{F}/F) }_{ {\rm Gal}(\overline{F}/F_n) }(\omega_{nf}^h) ) \otimes\mu_{\delta}$$
 is then uniquely determined by the two conditions 
 $$\det  \ind (\omega_{nf}^h) =  \omega^h_f\hskip10pt \text{and} \hskip10pt\ind (\omega_{nf}^h) |_{\cI} = \omega_{nf}^{h} \oplus \omega_{nf}^{qh} \oplus ...\oplus \omega_{nf}^{q^{n-1}h} .$$ Let $k/\bbF_q$ be a finite extension. Every absolutely irreducible smooth $k$-representation of ${\rm Gal}(\overline{F}/F) $ of dimension $n$ is isomorphic to 
$\ind (\omega_{nf}^h) \otimes \mu_{\lambda}$ for a $q$-primitive 
$1\leq h \leq q^n-2$ and a scalar $\lambda\in\overline{\bbF}_q^{\times}$ such that $\lambda^n\in k^{\times}$ and one has 
$$\ind (\omega_{nf}^h) \otimes \mu_{\lambda}\simeq \ind (\omega_{nf}^{\tilde{h}}) \otimes \mu_{\tilde{\lambda}}$$
if and only if  ${\rm Gal}(F_n/F).\omega_{nf}^h={\rm Gal}(F_n/F).\omega_{nf}^{\tilde{h}}$
and $\lambda^n=\tilde{\lambda}^n$.

\vskip5pt 

The theory of Lubin-Tate $(\varphi,\Gamma)$-modules and their relation to Galois representations is developed in \cite{KR09} and \cite{Sch17}. 
Let $F_{\phi}$ be a Lubin-Tate group for $\pi$, with Frobenius power series $\phi(t)\in o_F[[t]]$. The corresponding ring homomorphism $o_F\rightarrow \End(F_{\phi})$ is denoted by $a\mapsto [a](t)= at +...$. In particular, $[\pi](t)=\phi(t)$.
Let $F_{\infty}/F$ be the extension generated by all torsion points of $F_{\phi}$ and let $$H_F:={\rm Gal}(\overline{F}/F_{\infty})\hskip10pt \text{and} \hskip10pt
\Gamma:={\rm Gal}(\overline{F}/F)/ H_F = {\rm Gal}(F_{\infty}/F).$$ 

\vskip5pt
Let $z=(z_j)_{j\geq 0}$ be a $o_F$-generator of the Tate module of $F_{\phi}$. In particular, for $j\geq 0$
$$z_j=[\pi](z_{j+1})\equiv z_{j+1}^q\mod \pi$$ and $N_{F(z_1)/F}(-z_1)=\pi$. This implies
$$z_{j+1}^q=z_j(1+O(\pi)) \text{~for~} j\geq 1 \text{~~and~~}
z_1^{q-1} = -\pi (1+O(z_1)).$$

The Galois action on the generator $z$ is given by a character 
$\chi_L: {\rm Gal}(\overline{F}/F) \rightarrow o_F^{\times}$, which is surjective and has kernel $H_F$. One has $\omega_f \equiv \chi_L {\rm ~mod~ } \pi$.

\vskip5pt 

 The power series ring $o_F[[t]]$ has a Frobenius endomorphism and a $\Gamma$-action via
$\varphi(f)(t)=f( [\pi](t) )$ and $(\gamma f)(t)=f([\chi_L(\gamma)](t))$ for $f(t)\in o_F[[t]]$. Via reduction mod $\pi$, these actions induce a Frobenius action and a $\Gamma$-action on $\bbF_q[[t]]$ and its quotient field $\bbF_q((t))$.
This allows to introduce an abelian tensor category of \'etale Lubin-Tate $(\varphi,\Gamma)$-modules over $\bbF_q((t))$. It turns out to be canonically 
equivalent to the category of continuous finite-dimensional $\bbF_q$-representations of ${\rm Gal}(\overline{F}/F)$, cf. \cite[1.6]{KR09}, \cite[3.2.7]{Sch17}. 

\vskip5pt

To explain the functor from $(\varphi,\Gamma)$-modules to Galois representations, we denote by $\bbC_p$ the completion of an algebraic closure of $\bbQ_p$ and choose an embedding $\overline{F}\subseteq \bbC_p$. Recall
that the tilt $\bbC_p^{\flat}$ of the perfectoid field $\bbC_p$ is an algebraically closed and perfect complete non-archimedean field of characteristic $p$. Its valuation ring $o_{\bbC^{\flat}_p}$ is given by the projective limit
$\varprojlim_{x\mapsto x^q} o_{\bbC_p}/\pi o_{\bbC_p}$ and its residue field is $\overline{\bbF}_q$. 
There is a unique multiplicative section $$s: \overline{\bbF}_q\longrightarrow o_{\bbC^{\flat}_p}, a \mapsto (\tau(a){\rm ~mod~} \pi, \tau(a^{q^{-1}}){\rm ~mod~} \pi, \tau(a^{q^{-2}}){\rm ~mod~} \pi,...)$$ 
where $\tau$ denotes the Teichm\"uller map $\overline{\bbF}_{q}\rightarrow o_{\bbC_p}$. There is an inclusion $$\bbF_q((t))\stackrel{\subset}{\longrightarrow}  \bbC_p^{\flat}, ~ t\mapsto  (...,z_j {\rm ~mod}~ \pi,...)$$ and one has $\bbC_p^{\flat}= o_{\bbC^{\flat}_p} [1/t]$. The field $\bbC_p^{\flat}$ is endowed with a continuous action of ${\rm Gal}(\overline{F}/F)$ and 
a Frobenius $\varphi_q$, which raises any element to its $q$-th power.  We let $\bbF_q((t))^{\rm sep}$ denote the separable algebraic closure of $\bbF_q((t))$ inside $\bbC_p^{\flat}$. The field $\bbF_q((t))$ and its separable closure $\bbF_q((t))^{\rm sep}$ inherit the Frobenius action and the commuting ${\rm Gal}(\overline{F}/F)$-action from $\bbC_p^{\flat}$ and there is an
isomorphism $$H_F\stackrel{\simeq}\longrightarrow {\rm Gal}(\bbF_q((t))^{\rm sep}/\bbF_q((t))).$$

The functor $\sV$ from  $(\varphi,\Gamma)$-modules to Galois representations is then given by $$D \rightsquigarrow \sV(D):=(\bbF_q((t))^{\rm sep}\otimes _{\bbF_q((t))} D)^{\varphi=1}$$ 
where ${\rm Gal}(\overline{F}/F)$ acts diagonally (and via its projection to $\Gamma$ on the second factor).

\vskip5pt 

Now suppose that $k/\bbF_q$ is a finite extension. One can consider a $k$-representation of ${\rm Gal}(\overline{F}/F)$ as an $\bbF_q$-representation with a $k$-linear structure. Similarly, one may introduce $(\varphi,\Gamma)$-modules over $k((t))=k\otimes_{\bbF_q} \bbF_q((t))$, where $k$ has the trivial Frobenius and $\Gamma$-action. The functor $\sV$ then restricts to an equivalence of categories between  
\'etale $(\varphi,\Gamma)$-modules over $k((t))$ and continuous finite-dimensional $k$-representations of ${\rm Gal}(\overline{F}/F)$.

\section{The main result}

We fix once and for all an element $y\in \bbF_q((t))^{\rm sep}$ such that 

$$ y^{(q^n-1)/(q-1)}=t.$$

 For $g\in {\rm Gal}(\overline{F}/F)$, the power series $$f_g(t)=\chi_L(g)t/g(t)\in 1+ (\pi) [[t]]$$
depends only on the class of $g$ in $\Gamma$. The same is true for its mod $\pi$ reduction $\overline{f}_g(t)=\omega_{f}(g)t/g(t).$ Note that the formula  $f^s_g(t)$ defines an element of $o_F[[t]]$ for any $s\in\bbZ_p$.

\begin{Lem}\label{lem-aux} One has $g(y)=y \omega_{nf}^q(g)\overline{f}_g^{-\frac{q-1}{q^n-1}}(t)$ in $\bbF_q((t))^{\rm sep}$ for all $g\in {\rm Gal}(\overline{F}/F_{n})$. 
\end{Lem}
\begin{proof}
This is a version of \cite[Lem. 2.1.3]{Be10}. Let $j\geq 1$ and choose
$\pi_{nf,j}\in o_{\bbC_p}$ such that $$\pi_{nf,j}^{\frac{q^n-1}{q-1}}=z_j.$$
We write $\pi_j$ for $\pi_{nf,j}$ in the following calculations. Let $g\in {\rm Gal}(\overline{F}/F_{n})$. 
Then $$(g(\pi_j)/\pi_j)^{\frac{q^n-1}{q-1}}=g(z_j)/z_j=\chi_L(g)f^{-1}_g(z_j)$$ and so the quotient of 
$g(\pi_j)/\pi_j$ by $f_g^{-\frac{q-1}{q^n-1}}(z_j)$ is a certain ${\frac{q^n-1}{q-1}}$-th root
of $\chi_L(g)$. Since exponentiation with ${\frac{q^n-1}{q-1}}\in\bbZ_p^{\times}$ is surjective on the subgroup $1+(\pi) \subset o_F^{\times}$
we may write this root as $\tau(\omega_{nf,j}(g))$, with an element $\omega_{nf,j}(g)\in \bbF_{q^n}^\times$, and arrive at 
$$ g(\pi_j)/\pi_j=\tau(\omega_{nf,j}(g))f_g^{-\frac{q-1}{q^n-1}}(z_j).$$

The map $g\mapsto \omega_{nf,j}(g)$ is a character of the group
 ${\rm Gal}(\overline{F}/F_{n})$, since  $$\omega_{nf,j}(g) \equiv g(\pi_j)/\pi_j \mod \mathfrak{m}_{\bbC_p}$$ in the field $\overline{\bbF}_q= o_{\bbC_p}/ \mathfrak{m}_{\bbC_p}$ and this element is fixed by ${\rm Gal}(\overline{F}/F_{n})$.
 Moreover, this character does not depend on the choice of $\pi_j$: a different choice $\pi'_j$ differs from $\pi_j$ by a ${\frac{q^n-1}{q-1}}$-th root of unity, i.e. by an element of $F_{n}$.
 Hence $g(\pi'_j)/\pi'_j=g(\pi_j)/\pi_j$. By this independence, we see (using the element $\pi_{j+1}^q$ as an alternative choice for $\pi_j$) that $$\omega_{nf,j+1}^q=\omega_{nf,j}~ \text{for} ~ j\geq 1.$$
 Moreover, $\pi_{nf,1}^{q^n-1}=z_1^{q-1} = -\pi (1+O(z_1))$ and so $(\pi_{nf}/ \pi_{nf,1})^{q^n-1} \equiv 1 \mod \mathfrak{m}_{\bbC_p}$. The quotient $\pi_{nf}/ \pi_{nf,1} \mod \mathfrak{m}_{\bbC_p}$ is therefore fixed by  ${\rm Gal}(\overline{F}/F_{n})$, in other words $$g(\pi_{nf,1})/\pi_{nf,1} \equiv g(\pi_{nf})/\pi_{nf} \mod \mathfrak{m}_{\bbC_p}$$ for all  $g \in {\rm Gal}(\overline{F}/F_{n})$ and so
 $$  \omega_{nf,1}=\omega_{nf}.$$
 
Now recall that there is an isomorphism $ \varprojlim_{x\mapsto x^q} o_{\bbC_p}\simeq o_{\bbC_p^{\flat}}$ of multiplicative monoids 
given by reduction modulo $\pi$. We use the notation $u=(u^{(j)})$ for elements in the projective limit $\varprojlim_{x\mapsto x^q} o_{\bbC_p}$. The element 
$y\in o_{\bbC_p^{\flat}}$ is given by  $(...,\pi_j {\rm ~mod}~ \pi o_{\bbC_p},...)$. Its preimage
$(y^{(j)})$ under the above isomorphism is therefore given by 
$y^{(j)}=\lim_{m\rightarrow \infty} \pi_{j+m}^{q^m}$. 
By the same argument, the preimage of the element $\overline{f}_g^{-\frac{q-1}{q^n-1}}(t)$ has coordinates

$$\overline{f}_g^{-\frac{q-1}{q^n-1}}(t)^{(j)}=\lim_{m\rightarrow \infty} (f_g^{-\frac{q-1}{q^n-1}}(z_{j+m}))^{q^m}.$$

The composite map $s: \overline{\bbF}_{q}\rightarrow  o_{\bbC_p^{\flat}} \simeq \varprojlim_{x\mapsto x^q} o_{\bbC_p}$, which we also denote by $s$, is given by
$a \mapsto (\tau(a), \tau(a^{q^{-1}}), \tau(a^{q^{-2}}),...)$.
Since $$s(\omega_{nf}(g)^q)^{(j)}=\tau(\omega_{nf}(g)^{q^{-j+1}})=\tau(\omega_{nf,j}(g)),$$ we may put everything together and obtain
$$ \frac{g(y^{(j)})}{y^{(j)}}=\lim_{m\rightarrow \infty} (\frac{g(\pi_{j+m})}{\pi_{j+m}})^{q^m}=\tau( \omega_{nf,j}(g))\lim_{m\rightarrow \infty} (f_g^{-\frac{q-1}{q^n-1}}(z_{j+m}))^{q^m}=s(\omega_{nf}(g)^q)^{(j)}
\overline{f}_g^{-\frac{q-1}{q^n-1}}(t)^{(j)}.$$ Reducing this equation modulo $\pi$ yields the assertion of the lemma. 
\end{proof}

We now consider the
 $(\varphi,\Gamma)$-modules associated to the
irreducible Galois representations of the form $\ind(\omega_{nf}^h)$. 

\begin{Th} The Lubin-Tate $(\varphi,\Gamma)$-module associated to an irreducible Galois representation of the form $\ind(\omega_{nf}^h)$
is defined over the ring $\bbF_q((t))$ and admits a basis $e_0,e_1,...,e_{n-1}$ in which 
$$\gamma(e_j)=\overline{f}_{\gamma}(t)^{hq^j(q-1)/(q^n-1)}e_j$$ for all $\gamma\in\Gamma$ and $\varphi(e_j)=e_{j+1}$ and $\varphi(e_{n-1})=(-1)^{n-1}t^{-h(q-1)}e_0$.
\end{Th}
\begin{proof}
Let $D$ be the $(\varphi,\Gamma)$-module described in the statement and let $W=\sV(D)$. With $x=t^h e_0 \wedge ... \wedge e_{n-1}$, one has
$$ \varphi(x)= \varphi(t)^h (-1)^{n-1}t^{-h(q-1)}  e_1 \wedge ... \wedge e_{n-1} \wedge e_0 = t^{qh-h(q-1)} e_0 \wedge ... \wedge e_{n-1} = x.$$
Moreover, 
$$ \gamma(t)^h \prod_{ j=0}^{n-1} \overline{f}_{\gamma}^{h q^j (q-1)/(q^n-1)}(t) = 
( \omega_f(\gamma) t / \overline{f}_{\gamma}(t) ) ^h \overline{f}_{\gamma}^{h(q-1)/(q^n-1) \sum_{ j=0}^{n-1} q^j}  = \omega_f(\gamma)^h t^h$$
which implies $\gamma(x)= \omega_f(\gamma)^h x$ for all $\gamma\in\Gamma$. So $\det W = \omega_f ^h$. Put $k:=\bbF_{q^n}$ as a coefficient field, i.e. endowed with the trivial Frobenius action. To complete the proof, it remains to check that the restriction of $k\otimes_{\bbF_q} W$ to the inertia subgroup $\cI$ is given by 
$\omega_{nf}^{h} \oplus \omega_{nf}^{qh}\oplus ...\oplus \omega_{nf}^{q^{n-1}h}$. There is a ring isomorphism  
$$k \otimes_{\bbF_q}\bbF_q((t))^{\rm sep}\stackrel {\simeq}{\longrightarrow}\prod_{j=0}^{n-1}\bbF_q((t))^{\rm sep},\; x\otimes z \mapsto (\varphi_q^j (x)z)$$ where $\varphi_q$ is the $q$-Frobenius on $k$.
The induced Frobenius and ${\rm Gal}(\overline{F}/F_{n})$-action on $\prod_{j=0}^{n-1} \bbF_q((t))^{\rm sep} $ are given as 
$$
\begin{array}{ccl}
 \varphi ((x_0,...,x_{n-1})) & =&  (\varphi_q (x_{n-1} ), \varphi_q(x_0),...,\varphi_q (x_{n-2})) \\
&&\\
g ((x_0,...,x_{n-1}))&=  &(g (x_0),...,g( x_{n-1} )). 
\end{array}
$$
Choose $\alpha\in\overline{\bbF}_q \subset \bbF_q((t))^{\rm sep} $ such that $\alpha^{q^n-1}=(-1)^{n-1}$ and define the elements 
$$
\begin{array}{ccl}
 v_0 & =&  (\alpha y^h,0,...,0)e_0+ (0, \alpha^q y^{qh},0,...,0)e_1+...+  (0,...,0, \alpha^{q^{n-1}} y^{q^{n-1}h})e_{n-1}\\
&&\\
v_1 &=  & (0,\alpha y^h,0,...,0)e_0+ (0,0, \alpha^q y^{qh},0,...,0)e_1+...+  (\alpha^{q^{n-1}} y^{q^{n-1}h},0,...,0)e_{n-1} \\
\vdots &&\\
v_{n-1} &=  & (0,...0,\alpha y^h)e_0+ (\alpha^q y^{qh},0,...,0)e_1+...+  (0,...,\alpha^{q^{n-1}} y^{q^{n-1}h},0)e_{n-1} . 
\end{array}
$$
By definition of $D$, the vectors $e_i$ form a $\bbF_q((t))$-basis for $D$ and it follows easily 
that the vectors $v_i$ form a $k \otimes_{\bbF_q}  \bbF_q((t))^{\rm sep} $-basis
 for $k\otimes_{\bbF_q} ( \bbF_q((t))^{\rm sep} \otimes_{\bbF_q((t))}  D )$. Moreover, 
 
 $$ 
 \begin{array}{ccl}
 
 \varphi (  (0,...,0, \alpha^{q^{n-1}} y^{q^{n-1}h})e_{n-1}) &= & ( \alpha^{q^{n}} y^{q^{n}h},0,...,0)  \varphi ( e_{n-1}) \\
 &&\\
&=&  ( \alpha^{q^{n}} y^{q^{n}h},0,...,0)) (-1)^{n-1}t^{-h(q-1)} e_{0} =  (\alpha y^h,0,...,0)e_0
 \end{array}
 $$
 since  $\alpha^{q^{n}} = (-1)^{n-1} \alpha$ and $y^{q^n} t^{1-q}=y$. This means 
  $$
\begin{array}{ccl}
 \varphi (v_0) & =&  (0,\alpha^q y^{qh},0,...,0) \varphi (e_0)+ ( 0,0, \alpha^{q^2} y^{q^2h},0,...,0) \varphi (e_1)+...+ ( \alpha^{q^{n}} y^{q^{n}h},0,...,0) \varphi (e_{n-1}) \\
&&\\
&=  &   (0,\alpha^q y^{qh},0,...,0) e_1+ ( 0,0, \alpha^{q^2} y^{q^2h},0,...,0) e_2+...+  (\alpha y^h,0,...,0)e_0 \\
&=& v_0.
\end{array}
$$
Similarly, one shows $\varphi(v_j)=v_j$ for $j\geq 1$, so that $$v_0,...,v_{n-1} \in k \otimes_{\bbF_q} ( \bbF_q((t))^{\rm sep}  \otimes_{\bbF_q((t))}  D )^{\varphi=1}= k \otimes_{\bbF_q} \sV(D)= k \otimes_{\bbF_q}  W.$$
 Now if $g\in{\rm Gal}(\overline{F}/F_{n})$, then $g(y)=y \omega_{nf}^q(g)c_g$ with $c_g:=\overline{f}_g^{-\frac{q-1}{q^n-1}}(t)$ by lemma \ref{lem-aux} and
 $ g(e_j)=c_g ^{-q^j h}e_j$ by definition of $D$. Hence
 $$ g(y)^{q^jh}g(e_{j})= (y \omega_{nf}^q(g) )^{q^jh}e_j.$$ 
 If $g\in\cI$, then $g(\alpha)=\alpha$ and then altogether
 
 $$
\begin{array}{ccl}
 g (v_0) & =&  (\alpha g(y)^{h},0,...) g (e_0)+ ( 0,\alpha^{q} g(y)^{qh},0,...) g (e_1)+...+ (0,...,\alpha^{q^{n-1}} g(y)^{q^{n-1}h}) g (e_{n-1}) \\
&&\\
&=  & \omega_{nf}^{qh}(g)\cdot (  (\alpha y^{h},0,...) e_0 + ( 0,\alpha^{q} y^{qh},0,...) e_1+...+ (0,...,\alpha^{q^{n-1}} y ^{q^{n-1}h}) e_{n-1})   \\
&&\\
&=&  \omega_{nf}^{qh}(g)\cdot v_0,
\end{array}
$$
where $\cdot$ refers to the left $k$-structure of $\prod_{j=0}^{n-1} \bbF_q((t))^{\rm sep}$.
Similarly, one shows $g(v_j)=  \omega_{nf}^{q^{1-j}h}(g)v_j$ for all $j\geq 1$ and $g\in\cI$. Since $\omega_{nf}^{q^n}= \omega_{nf}$ and hence $\omega_{nf}^{q^{1-j}h}= \omega_{nf}^{q^{n+1-j}h}$, this proves that 
the restriction of $k \otimes_{\bbF_q} W$ to $\cI$ is given by the sum of the characters
$\omega_{nf}^{h}, \omega_{nf}^{qh},...,\omega_{nf}^{q^{n-1}h}$. 
\end{proof}

As explained above, one may pass from irreducible representations of the form $\ind(\omega_{nf}^h)$ to the general case by twisting with characters. 
Note that any character ${\rm Gal}(\overline{F}/F)\rightarrow \overline{\bbF}^{\times}_q$ can be written in the form 
 $\omega_f^s \mu_{\lambda}$, for $1\leq s \leq q-1$ and $\lambda\in\overline{\bbF}^{\times}_q $.
\begin{Lem} Let $k/\bbF_q$ be a finite extension. The $(\varphi,\Gamma)$-module associated a Galois character of the form 
$\omega_f^s \mu_{\lambda}$ with $\lambda\in k^{\times}$ admits a basis $e$ such that $\varphi(e) = \lambda \cdot e$ and $\gamma (e) = \omega_f^s (\gamma)\cdot e$ for all $\gamma\in\Gamma$.
\end{Lem}
\begin{proof}
Since the functor $\sV$ preserves the tensor product, we may discuss the two characters  $\omega_f^s$ and $\mu_{\lambda}$ separately. 
For the twists by a character of $\Gamma$, such as $\omega_f^s$, see \cite[Remark 4.6]{SV16}.
So let $V=\mu_{\lambda}=k$ and let 
$$D(V)=(\bbF_q((t))^{\rm sep}\otimes _{\bbF_q} V)^{H_F}  $$ be the associated $(\varphi,\Gamma)$-module. 
It is instructive to check the case $k=\bbF_{q}$ first. Here, we choose $\beta\in \overline{\bbF}_{q}$ with $\beta^{q-1}=\lambda$ and 
put $e= \beta \otimes 1\in  \bbF_q((t))^{\rm sep} \otimes_{\bbF_q}  V$. Since $\beta \neq 0$, we have $e\neq 0$. Moreover, $\cI$ acts trivial on $e$ and for $\varphi\in {\rm Gal}(\overline{F}/F)$
we find 
$$ \varphi ( e ) = \varphi ( \beta ) \otimes \varphi ( 1 ) = \beta^q \otimes \lambda^{-1} = \beta \lambda \otimes \lambda^{-1} = \beta \otimes 1 = e.$$
Hence $e$ is indeed ${\rm Gal}(\overline{F}/F)$-invariant. Moreover, if $\phi$ denotes the Frobenius endomorphism on $D(V)$ we have 
$$\phi (e ) = \phi ( \beta ) \otimes 1 = \beta^q \otimes 1 =  \lambda \beta \otimes 1 = \lambda e.$$
Now suppose that $k=\bbF_{q^n}$ for some $n$ and $\lambda\in k^{\times}$. We use the ring isomorphism
$$k \otimes_{\bbF_q}\bbF_q((t))^{\rm sep}\stackrel {\simeq}{\longrightarrow}\prod_{j=0}^{n-1}\bbF_q((t))^{\rm sep},\; x\otimes z \mapsto (\varphi_q^j (x)z)$$ where $\varphi_q$ is the $q$-Frobenius on $k$. It is ${\rm Gal}(\overline{F}/F_n)$-equivariant, where the Galois action on the right-hand side is componentwise (see proof of the above theorem). By the normal basis theorem, there is $x\in k^{\times}$ such that its conjugates $\varphi_q^j (x)$ are linearly independent over ${\bbF_q}$.
The $j$-th copy $\bbF_q((t))^{\rm sep}$ in the above product has therefore a $\bbF_q((t))^{\rm sep}$-basis element $e_j := \varphi_q^j (x)\in k=V$ on which $\cI$ acts trivial and on which the element $\varphi^n\in {\rm Gal}(\overline{F}/F_n)$ acts by $\lambda^{-n}$. Choose $\beta\in \overline{\bbF}_{q}$ such that $\beta^{q^n-1}=\lambda^{n}$ and put 
$v_j = \beta e_j $. Then $\cI$ obviously acts trivial on $v_j$ and the same holds for $\varphi^n$, since 
$$\varphi^n(v_j) = \varphi^n (\beta) \varphi^n (e_j ) =  \beta^{q^n} \lambda^{-n} e_j  = \beta \lambda^{n}  \lambda^{-n} e_j =v_j.$$
Hence, $\cI$ and $\varphi^n$ act trivial on $(v_j)\in \prod_{j=0}^{n-1}\bbF_q((t))^{\rm sep}$ and then also on its preimage $v=x \otimes \beta\in  k \otimes_{\bbF_q}\bbF_q((t))^{\rm sep}$. 
Note that $v\neq 0$ since $x,\beta \neq 0$. Write $N= \prod_{j=0}^{n-1} \varphi^j$ and $e= N(v)$. Then $e$ is fixed by $\cI$ (since $\cI$ is normalized by the $\varphi^j$) and is fixed by $\varphi$ by construction. Hence, $e$ is ${\rm Gal}(\overline{F}/F)$-invariant. Note that $e\neq 0$, since $e=N(x) \otimes N(\beta)$ and $N(x),N(\beta) \neq 0$ and so $e$ is indeed a basis element of $D(V)$ on which $\Gamma$ acts trivial. Finally, write $e=\sum_{j=0}^{n-1}  \varphi^j_q(x)\otimes z_j $
with $z_j \in \bbF_q((t))^{\rm sep}$. The Frobenius endomorphism $\phi$ on $D(V)$ satisfies
$$\phi ( e ) = \sum_j \varphi^j_q(x)\otimes \varphi (z_j)  = \varphi ( \sum_j  \varphi^{-1} (\varphi^j_q(x)) \otimes z_j )= 
 \varphi ( \sum_j  \lambda\varphi^j_q(x)\otimes z_j ) = \lambda \varphi (e) =\lambda e. $$ 
 \end{proof}

\begin{Cor} 
Let $k/\bbF_q$ be a finite extension. The $(\varphi,\Gamma)$-module associated to an irreducible Galois representation of the form 
$( \ind(\omega_{nf}^h)) \otimes \omega_f^s \mu_{\lambda}$, for $1\leq s \leq q-1$ and $\lambda^n \in k^{\times}$,
is defined over the ring $k((t))$ and admits a basis $e_0,e_1,...,e_{n-1}$ in which 
$$\gamma(e_j)=\omega_f(\gamma)^s \overline{f}_{\gamma}(t)^{hq^j(q-1)/(q^n-1)}e_j$$ for all $\gamma\in\Gamma$ and 
$\varphi(e_j)=\lambda e_{j+1}$ and $\varphi(e_{n-1})=(-1)^{n-1}t^{-h(q-1)}\lambda e_0$.
\end{Cor}
\begin{proof}
This follows from the preceding lemma and the theorem.
\end{proof}

 Since $\omega_{nf}^{\frac{q^n-1}{q-1}}=\omega_f$, every irreducible representation of ${\rm Gal}(\overline{F}/F) $ of dimension $n$ is therefore isomorphic to 
$\ind (\omega_{nf}^h) \otimes \omega_f^s \mu_{\lambda}$ for $1\leq s \leq q-1$, a scalar $\lambda\in \overline{\bbF}^{\times}_q$ and 
a $q$-primitive $1 \leq h\leq \frac{q^n-1}{q-1} -1.$

\vskip10pt 

\noindent {\small Cédric Pépin, LAGA, Université Paris 13, 99 avenue Jean-Baptiste Clément, 93 430 Villetaneuse, France \newline {\it E-mail address: \url{cpepin@math.univ-paris13.fr}} }

\vskip10pt

\noindent {\small Tobias Schmidt, Univ Rennes, CNRS, IRMAR - UMR 6625, F-35000 Rennes, France \newline {\it E-mail address: \url{tobias.schmidt@univ-rennes1.fr}} }

\end{document}